\documentclass{amsart}

\usepackage{ifthen}

\newtheorem{anyprop}{Anyprop}[section]

\newtheorem{theorem}[anyprop]{Theorem}
\newtheorem{lemma}[anyprop]{Lemma}
\newtheorem{proposition}[anyprop]{Proposition}
\newtheorem{corollary}[anyprop]{Corollary}

\theoremstyle{definition}

\newtheorem{definition}[anyprop]{Definition}

\newtheorem{example}[anyprop]{Example}

\newtheorem{remark}[anyprop]{Remark}

\newcommand{\CC}{\mathbb{C}}

\newcommand  {\shL}     {\mathcal{L}}

\newcommand  {\shS}     {\mathcal{S}}

\newcommand  {\fom}     {\mathfrak{m}}

\newcommand  {\fon}     {\mathfrak{n}}

\newcommand  {\fop}     {\mathfrak{p}}

\renewcommand  {\ker }  {\operatorname{ker}}

\newcommand  {\N}       {\operatorname{N}}

\renewcommand{\O}       {\mathcal{O}}

\newcommand  {\Pic}     {\operatorname{Pic}}

\newcommand  {\rk}    {\operatorname{rk}}

\newcommand  {\Spec}    {\operatorname{Spec}}

\newcommand  {\Sym}     {\operatorname{Sym}}

\newcommand  {\Syz}     {\operatorname{Syz}}

\newcommand{\comdots}{ , \ldots , }

\newcommand{\plusdots}{ + \ldots + }

\theoremstyle{remark}

\numberwithin{equation}{section}

\usepackage{amscd}
\usepackage{amssymb}
\input xy
\xyoption{all}

\newcommand{\mathdisp}[2]{\[ #1 \, #2 \]}

\newcommand{\mathl}[2]{$#1$#2}

\newcommand{\betonung}[1]{\emph{#1}}

\newcommand{\forca}{{B}}

\newcommand{\repmatrix}{{\alpha}}

\newcommand{\repmentry}{\alpha}

\newcommand{\classnen }{\beta}

\begin{document}
\title[Some remarks on the affineness of ${\mathbb A}^1$-bundles]
{Some remarks on the affineness of ${\mathbb A}^1$-bundles}

\author[Holger Brenner]{Holger Brenner}
\address{Universit\"at Osnabr\"uck, Fachbereich 6: Mathematik/Informatik,
Albrechtstr. 28a,
49069 Osnabr\"uck, Germany}
\email{hbrenner@uni-osnabrueck.de}


\subjclass{}


\dedicatory{Meinem akademischen Lehrer, Herrn Prof. Dr. Uwe Storch zu seinem 70. Geburtstag gewidmet.}

\begin{abstract}
We study how forcing algebras give rise to ${\mathbb A}^1$-bundles and ${\mathbb A}^1$-torsors and how they are related to ${\mathbb A}^1$-patches. In particular we discuss the affineness of torsors and how algebraic properties of ${\mathbb A}^1$-patches can be deduced from this property.
\end{abstract}

\maketitle

\noindent Mathematical Subject Classification (2010): 13 A 35, 14 R 10, 14 R 25

\smallskip

\noindent Keywords: Forcing algebra, torsor, affine scheme, cancellation problem, ${\mathbb A}^1$-bundle, ${\mathbb A}^1$-fibration

\section*{Introduction}

In many different situations one is interested in schemes which allow some kind of an ${\mathbb A}^1$-fibration: group actions of the additive or the multiplicative group, locally nilpotent derivations, line bundles and their torsors, cancellation problems.

In this article, which is in part a survey article, we want to relate these notions with the concept of a forcing algebra, a concept introduced by Mel Hochster \cite{hochstersolid}. A forcing algebra has the form
\[ B= R[T_1, \ldots ,T_n]/(f_1T_1 + \ldots + f_nT_n +f) \, ,  \]
where $R$ is a commutative ring and $f,f_1, \ldots ,f_n \in R$. They are mainly used to understand closure operations for ideals: whether the element $f$ belongs to some closure (like the radical, integral closure, tight closure) of the ideal $(f_1, \ldots ,f_n)$  can often be characterized by properties of the $R$-algebra $B$. For a point $P \in \Spec \, R$, the fiber over $P$ is the solution set of the inhomogeneous linear equation
\[ f_1(P)T_1 + \ldots + f_n(P)T_n +f(P) = 0\]
over $\kappa(P)$. For $n=2$, this solution set is generically an affine line, but it can degenerate to an affine plane or the empty set. 

In this paper we want to focus on the ${\mathbb A}^1$-properties of the corresponding schemes $\Spec B$ and open subsets whereof, especially for two ideal generators ($n=2$). Over the open subset $U=D(f_1, \ldots, f_n)$ the structure of $(\Spec B)|_U$ is quite simple, it is an affine-linear bundle in general, and the sheaf of syzygies $\shS= \Syz(f_1, \ldots , f_n)$, which is locally free on $U$, acts on it by translation. These affine-linear bundles or $\shS$-torsors are classified by the first cohomology group of the syzygy sheaf. When starting with such a cohomology class and its corresponding torsor, a forcing algebra yields a natural affine completion of the torsor. For $n=2$, we get $\shS \cong \O$ on $U$ and so we are dealing with ${\mathbb A}^1$-torsors.

The relation between torsors and forcing algebras is often helpful in understanding global properties of torsors. Here we are in particular interested in the quesion whether the torsor is an affine scheme and whether its ring of global sections is finitely generated.

Let me give a quick overview on the organization of this paper. In Section \ref{forcingsection} we introduce forcing algebra given by a (finitely generated) submodule and an element in some (finitely generated) module. In Sections \ref{torsorsection} and \ref{fibrationsection} we deal with torsors of a vector bundle, these are affine-linear bundles and they are classified by the first cohomology group of the vector bundle which acts on them. In the easiest case, the acting vector bundle is just the affine line, giving rise to ${\mathbb A}^1$-torsors. A forcing algebra gives rise to such a torsor over a suitable open subset of $\Spec R$.

Section \ref{patchsection} deals with ${\mathbb A}^1$-patches, a concept introduced by Dutta, Gupta and Onoda \cite{duttaguptaonoda} and presented at the conference in Bangalore, with the focus on algebraic properties of these rings like flatness and finite generation. We show that such ${\mathbb A}^1$-patches arise in a natural way as rings of global sections of suitable ${\mathbb A}^1$-torsors. The main point is to understand whether these torsors are affine schemes or not. For a two-dimensional base ring, this question is surprisingly related to the theory of tight closure, as explained in Section \ref{backgroundsection}. With this background at hand we can answer some of their questions in Section \ref{affinesection}.

In Section \ref{cechsection} we describe how the torsor defined by a ${\rm \check{C}}$ech cohomology class of the structure sheaf on a quasiaffine scheme can be realized by a forcing algebra defined by a system of forcing equations. In Section \ref{cancellationsection} we recall a result of D.R. Finston and S. Maubach \cite{finstonmaubach} on how torsors and their affineness properties are related to the cancellation problem, i.\,e.\  the problem whether $X \times {\mathbb A}^1 \cong {\mathbb A}^{n+1}$ implies $X  \cong {\mathbb A}^{n}$. Based on these two sections we establish in Section \ref{duboulozsection}, which started from discussions with A. Dubouloz in Bangalore and which is related to an example of J. Winkelmann \cite{winkelmannfreeholomorphic}, a class of examples of torsors which are counter-candidates to the cancellation problem. 

In Section \ref{derivationsection} we show that a forcing algebra for two ideal generators defining an ideal of depth two carries a locally nilpotent derivation and in Section \ref{singularsection} we discuss some smoothness properties.

This paper arose from the conference held in December 2010 in Bangalore in honor of R.V. Gurjar, Balwant Singh and Uwe Storch. During this conference I had the opportunity to discuss with many people on topics related to this paper, in particular with S.M. Bhatwadekar, A. Dubouloz, A. Dutta, R.V. Gurjar, H. Flenner, Neena Gupta, K. Masuda, M. Miyanishi, U. Storch. I thank them for these discussions and the organizers of the conference,  H. Flenner,  M. Miyanishi, D. Patil, R. A. Rao, L. G. Roberts  and J. K. Verma. I thank A. Dutta and N. Gupta for their hospitality during my stay in Kolkata in July 2011. I thank D. Gomez Ramirez for discussions on global properties of forcing algebras and N. Gupta and A. St\"abler for many corrections and suggestions.

I dedicate this paper to my academic teacher, Prof. Dr. Uwe Storch. I am grateful for the mathematical education I enjoyed in Bochum, for the mathematical spirit he represents and for his advice.

\section{Forcing algebras}
\label{forcingsection}

We start with the definition of a forcing algebra, for related material see \cite{brennerforcingsyzygy} and \cite{brennerkolkata}.

\begin{definition}
Let $R$ be a commutative ring and let $f_1, \ldots ,f_n, f \in R $ be elements. The $R$-algebra
\[B= R[T_1, \ldots ,T_n]/(f_1T_1 + \ldots + f_nT_n +f) \]
is called the \emph{forcing algebra} for the data $f_1, \ldots, f_n,f$.
\end{definition}

Forcing algebras occured first in the work of Hochster in relation to solid closure \cite{hochstersolid}. The universal property of the forcing algebra is that $f \in IB$, where  $I=(f_1, \ldots ,f_n)$, and that for every $R$-algebra $C$ with $f \in IC$ there exists a (non-unique) ring homomorphism $B \rightarrow C$. The question whether $f$ belongs to a certain ideal closure of $I$ (like the radical, integral closure, Frobenus closure, tight closure, plus closure, solid closure) can often be translated into certain properties of the corresponding forcing algebra. Typical examples are the following relations.
\begin{enumerate}
\item
\[ f \in I \text{ if and only if } \varphi: \Spec B \longrightarrow \Spec R \text{ has a section}. \]
Here it doesn't matter whether the section is a ring homomorphism or an $R$-module section (i.\,e.\ $R$ is a direct summand of $B$).

\item
The containment inside the radical can be expressed by
\[ {f \in {\rm rad }\, I  \text{ if and only if } \varphi \text{ is surjective}} \, {.}\]

\item
The containment inside the integral closure $\bar{I}$ can be characterised by
\[ {f \in \overline{ I}\, \,  \text{ if and only if } \varphi \text{ is a universal submersion}}  {.}\]
\end{enumerate}

Locally, $\Spec B$ has an easy description. For every point $P\ \in X=\Spec R$ the fiber over $\kappa (P)$ is just the solution set to the inhomogeneous linear equation
\[f_1(P) t_1 + \ldots + f_n(P)t_n + f (P) = 0 \, ,\]
which is empty or an $n-1$ or $n$-dimensional affine space. Over $D(f_i)$, $i=1 , \ldots , n$, we can write
\[B_{f_i} \cong R_{f_i}[T_1 , \ldots , T_{i-1}, T_{i+1}, \ldots , T_n] \, ,\]
so this is locally an affine space of dimension $n-1$ over $D(I)$. For a point $P \in V(I)$, all $f_i(P)$ vanish and the fiber over $P$ is either empty (if $f(P) \neq 0$) or has dimension $n$.

To study also closure operations on submodules, the following generalization of forcing algebras is helpful. Let
\[B = R[T_1, \ldots, T_n]/(\repmatrix  T-s) , \,\]
where $\repmatrix =(\repmentry_{ij})$ is an $m \times n$-Matrix over $R$, $T=(T_1, \ldots, T_n)$ and $s =(s_1, \ldots, s_m) \in R^m$. In more detail, such algebras are defined by a system of inhomogeneous forcing equations, namely
\[ \begin{matrix}
 \repmentry_{11}T_1 & +& \ldots &+  &  \repmentry_{1n}T_n &=& s_1  \cr
\repmentry_{21}T_1 & + & \ldots &+ &  \repmentry_{2n}T_n &=& s_2  \cr
 & & & & \cr
\repmentry_{m1}T_1 & +& \ldots &+ &  \repmentry_{mn}T_n &=& s_m     \,.
\end{matrix} \]

We explain how these algebras arise from a finitely generated submodule $N \subseteq M$ inside a finitely generated $R$-module $M$ and an element $s \in M$. Let $y_1
\comdots y_m$ be generators for $M$ and $x_1 \comdots x_n$ generators for $N$. This gives rise to a surjective homomorphism
$\psi: R^m \rightarrow M $, a submodule $N' = \varphi^{-1}(N)$ and a
morphism $R^n \rightarrow R^m$ which sends $e_i$ to a preimage
$x_i'$ of $x_i$. Altogether we get the commutative diagram with exact rows
$$ \xymatrix{ & & R^n &\! \! \!\stackrel{\repmatrix}{ \longrightarrow }\!\! \!  &
R^m &\! \! \! \longrightarrow \!\! \! & M/N & \!\! \! \longrightarrow \! \! \!& \! \! 0  \cr
& & \downarrow  & & \downarrow \psi & & \downarrow = &  &\cr
0\! \! &\! \! \!\!  \longrightarrow\! \! \! \! & N & \! \!\!\longrightarrow\!\! \! &M &\!\! \! \longrightarrow\! \! \! &M/N&\!\! \!\longrightarrow\! & 0  } $$
($\repmatrix$ is an $m \times n$-matrix). The element $s $ is represented by $(s_1 \comdots s_m) \in R^m$,
and $s $ belongs to $N$ if and only if the system
$$ \repmatrix \begin{pmatrix} t_1 \cr . \cr . \cr . \cr t_n \end{pmatrix}
  = \begin{pmatrix} s_1 \cr . \cr . \cr . \cr s_m \end{pmatrix}  $$
has a solution.

The formation of forcing algebras commutes with arbitrary base change $R \rightarrow R'$. Therefore for every point $P \in \Spec R$ the fiber ring $\forca \otimes_R \kappa(P)$ is the forcing algebra given by $$\repmatrix (P) T= s(P) \, ,$$ which is a system of inhomogeneous linear equations over the field $\kappa(P)$. Hence the fiber of $\Spec \forca \rightarrow \Spec R$ over $P$ is the \emph{solution set} to a system of linear inhomogeneous equations.

\section{Torsors}
\label{torsorsection}

A group scheme $G$ over a scheme $U$ is a scheme $G \rightarrow U$ together with a group operation $G \times_U G \rightarrow G$, an inverse mapping $G \rightarrow G$ and a (neutral) section $U \rightarrow G$ satisfying certain natural conditions. We consider a geometric vector bundle $V \rightarrow U$ together with its addition, its negation and its zero section as a group scheme.

\begin{definition}
Let $G \rightarrow U$ be a group scheme. A scheme $T$ over  $U$ together with a group action $G \times T \rightarrow T$  is called a $G$-\emph{torsor} (in the Zariski topology), if there exists an open covering of $U$ such that the action is for this covering isomorphic to the action of the group on itself.
\end{definition}

Torsors are also called principal homogeneous spaces or principal fiber bundles. The spectrum of a forcing algebra for $s=0$ is an affine commutative group scheme $G$ which acts on $T=\Spec B$. The restrictions of these objects to $U=D({\rm Ann} M/N))$ (to $U=D(I)$ in the case of an ideal $I \subseteq R$) is particularly important. The restriction $V=G |_{U}$ is a vector bundle (which we call the syzygy bundle $\Syz( \repmatrix)$ or $\Syz(f_1, \ldots , f_n)$ in the ideal case) and $T|_{U}$ is a $V$-torsor. Hence a forcing algebra yields a torsor over $U$ and the spectra of forcing algebras may be thought of as affine completions of torsors over quasiaffine schemes; see also Section \ref{cechsection}).

On a separated noetherian scheme the $V$-torsors are classified by $H^1(U, \shS)$, where $\shS$ is the locally free sheaf of sections in the vector bundle $V$, see \cite[Proposition 3.2]{brennerbarcelona}. This rests on the fact that the glueing data for a principal fiber bundle correspond to ${\rm \check{C}}$ech cohomology classes. We denote the torsor corresponding to a cohomology class $c$ by $T(c)$. The torsor is trivial if and only if it has a section if and only if it is isomorphic to the vector bundle $V$. A torsor is trivial over every affine open subscheme. Torsors have nice functorial properties, in particular, if $T(c) \rightarrow U$ is a torsor and $U' \rightarrow U$ is a morphism, then $T(c) \times_U U' \cong T(c')$, where $c'$ is the pull-back of $c$ (see  \cite[Proposition 3.3]{brennerbarcelona}).

It is easy to describe the cohomology class given by the torsor induced by forcing data.

\begin{proposition}
\label{restrictedforcing}
Let $N \subseteq M$ be finitely generated $R$-modules, let $s \in M$ and let $B=R[T_1, \ldots ,T_n]/( \repmatrix T-s)$ be a forcing algebra for these data. Let
$U$ be the open complement of the support of $M/N$ and let $\shS= \ker ( \O^n \stackrel{\repmatrix}{\rightarrow} \O ^m)$ on $U$. Then $\Spec B|_U \rightarrow U$ is the $\shS$-torsor on $U$ corresponding to $\delta(s) \in H^1(U,\shS)$.
\end{proposition}
\begin{proof}
Note that the map given by $\repmatrix$ is surjective on $U$ and therefore the syzygy sheaf $\shS=\Syz(\repmatrix )$ is a locally free sheaf on $U$. This gives a short exact sequence
\[ 0 \longrightarrow \shS \longrightarrow {\mathcal O}^n  \longrightarrow {\mathcal O}^m   \longrightarrow 0 \]
on $U$ and $s$ is represented by an element (also denoted by) $s \in \Gamma(U, {\mathcal O}^m )$. By the connecting homomorphism, $s$ defines a cohomology class $\delta(s) \in H^1(U,\shS)$. For the proof when $M=R$, $N=I$ an ideal (and hence $U=D(I)$), we refer to \cite[Proposition 3.5]{brennerbarcelona}. The same proof works also in the module case.
\end{proof}

There are situations where the fibers of $\Spec B$ over $\Spec R \setminus U$ are empty, e.\,g.\ if $IB$ is the unit ideal. In this case the spectrum of the forcing algebra itself equals the induced torsor, which is then an affine scheme. In general it is difficult to determine whether a torsor is an affine scheme (i.\,e.\ the spectrum of a ring).

\section{${\mathbb A}^1$-fibrations and $({\mathbb A}^1,+)$-actions}
\label{fibrationsection}

In this paper we will focus on forcing algebras such that their induced torsors are ${\mathbb A}^1$-torsors. In particular we look at forcing algebras \[B=R[T_1,T_2]/(f_1T_1+f_2T_2+f) \, , \] where $I=(f_1,f_2)$ is an ideal generated by a regular sequence of two elements $f_1$ and $f_2$ in a domain $R$. This means $R_{f_1} \cap R_{f_2} =R$ (inside the quotient field $Q(R)$) and $\O \cong \Syz(f_1,f_2) $ by sending $1 \mapsto (-f_2,f_1)$. So in this situation the syzygy bundle is just the structure sheaf and the torsor on $D(I)$ is given by the ${\rm \check{C}}$ech cohomology class $ \frac{f}{f_1f_2} \in H^1(U, \O)$. The action of the additive group scheme $G_a=({\mathbb A}^1_R,+) =\Spec R[W]$ (over $\Spec R$) is explicitly given by
\[ G_a \times_{\Spec R} \Spec B \longrightarrow  \Spec B,\,  (P, w, t_1,t_2) \longmapsto (P, t_1+f_2(P) w, t_2-f_1(P) w) \, .      \]
or, on the ring level, by
\[B \longrightarrow B[W] ,\, T_1 \longmapsto T_1+f_2 W,\, T_2 \longmapsto T_2 -f_1 W \, .\]

\begin{example}
Different forcing algebras may induce the same torsor on $D(I)$. A typical example is $R=K[x,y]$, ideals $(x^r,y^s)$ and $f =x^ay^b$. The forcing algebras
\[K[x,y][T_1,T_2]/(x^rT_1+y^sT_2 +x^ay^b)  \]
depend on $a,b,r,s$ (in fact this dependance is quite difficult, see Example \ref{torsoraffineplane} below), but the induced torsors on $D(I)=D(x,y)$ depend only on the cohomology classes $ \frac{x^a y^b}{x^ry^s} = \frac{1}{x^{r-a}y^{s-b} }  $ (if $r \geq a$ and $s \geq b$; else the class is $0$ anyway).
\end{example}

We recall some other related notions which apply to an ${\mathbb A}^1$-torsor.

\begin{definition}
An ${\mathbb A}^1$-fibration is an affine morphism of schemes $\varphi:Y \rightarrow X$ of finite type such that the fiber over every generic point of $X$ is an affine line.
\end{definition}

An ${\mathbb A}^1$-torsor $T \rightarrow U$ and a forcing algebra $\Spec B \rightarrow \Spec R$ as above (which induces such a torsor) are easy examples of ${\mathbb A}^1$-fibrations. Already these forcing algebras show that an ${\mathbb A}^1$-fibration might contain exceptional fibers.

\begin{example}
We look again at
\[B=K[x,y][T_1,T_2]/(x^rT_1+y^sT_2 +x^ay^b) \, . \]
All fibers over $D(x,y)$ are affine lines, since $(\Spec B)|_{D(x,y)}$ is just the corresponding ${\mathbb A}^1$-torsor $T$. The fiber over $(0,0)$ might have very different properties. If $a=b=0$, then the fiber over $(0,0)$ is empty, so in this case $T = \Spec B$. If $r=0$ or $s=0$, then we can eliminate $T_1$ (or $T_2$) to see that $\Spec B$ is just the affine line over ${\mathbb A}^2_K$. So let $r,s,a \geq 1$. Then the fiber over $(0,0)$ is an affine plane. Considered as a subset of $\Spec B$, it might contain singular points.
\end{example}

\section{${\mathbb A}^1$-patches}
\label{patchsection}

In \cite{duttaguptaonoda}, the authors introduced the notion of an ${\mathbb A}^1$-patch.

\begin{definition}
\label{a1patch}
Let $R$ be a noetherian domain and $A$ an $R$-domain. Then $A$ is called an ${\mathbb A}^1$-patch over $R$, if there exists a regular sequence $x,y \in R$ such that the following properties hold.
\begin{enumerate}
\item
$A_x = R_x[U]$.

\item
$A_y =R_y[V]$.

\item

$A=A_x \cap A_y$.
\end{enumerate}
\end{definition}
Here the intersection is of course to be taken inside the quotient field of $A$.

\begin{proposition}
\label{patchproperties}
Let $R$ be a noetherian domain and let $f_1, f_2 \in R$ be a regular sequence. Let $f \in R$ and consider the forcing algebra \[B=R[T_1,T_2]/(f_1T_1+f_2T_2+f) \, .\] Then the following hold.

(1) The open subset $D((f_1,f_2)B) \subseteq \Spec B$ is an ${\mathbb A}^1$-torsor over $D(f_1,f_2) \subseteq \Spec R$ corresponding to the ${\check{C}}$ech cohomology class $ \frac{f}{f_1f_2} \in H^1(D(f_1,f_2) , \O_R)$.

(2) The ring of global sections
\[A = \Gamma (D((f_1,f_2)B),\O_B)    \]
is an ${\mathbb A}^1$-patch.

(3) If $D((f_1,f_2)B) \subseteq \Spec B$ is an affine scheme, then $A$ is a finitely generated flat $R$-algebra.
\end{proposition}
\begin{proof}
(1) is a special case of Proposition \ref{restrictedforcing}\,(2). We have $D( (f_1,f_2)B) = D(f_1) \cup D(f_2)$ (in $\Spec B$) and we have $\Gamma (D(f_1), \O_B)=B_{f_1} \cong  R[T_2]$ and $\Gamma (D(f_2), \O_B)=B_{f_2} \cong  R[T_1]$. Since $B$ is a domain (because $f_1,f_2$ are regular elements), the ring of global sections of an open subset is just the intersection of the global rings for an open cover, hence \[\Gamma(D((f_1,f_2)B),\O_B) =  B_{f_1} \cap B_{f_2} \, .\]
Note also that $A_{f_i} \cong B_{f_i}$, since for a quasiaffine scheme $Z$ the canonical morpism $Z \rightarrow \Spec \Gamma(Z,\O_Z)$ is an open immersion.

(3) The ring of global sections of an affine open subscheme $U \subseteq \Spec B$ is always finitely generated over $B$ \cite[Satz 1.4.3]{brennerdissertation}. Therefore, in our situation, it is also finitely generated over $R$. The scheme morphism $D((f_1,f_2)B) \rightarrow D(f_1,f_2) \subseteq \Spec R$ is always flat, since the first morphism is locally an affine cylinder and the second is an open immersion. In the affine case we have $\Spec A \cong D((f_1,f_2)B) $, and so $A$ is flat over $R$.
\end{proof}

From this proposition we see that it is important to understand when a torsor is an affine scheme. Whenever this property is established, flatness and the property of being of finite type follow immediately. For the interplay of these notions see also \cite[Lemma 4.3]{duttaguptaonoda} and \cite[Corollary 3]{schenzelflatfinite}.
 
The following definition is a slight generalization of an ${\mathbb A}^1$-patch.

\begin{definition}
\label{a1patchgeneral}
Let $R$ be a commutative ring, $X= \Spec \,R$ and $A$ an $R$-algebra. Then $A$ is called an ${\mathbb A}^1$-patch over $R$, if there exists an ideal $I \subseteq R$ such that the following properties hold.
\begin{enumerate}
\item
$\Gamma(D(I), \O_X)=R$.

\item
$(\Spec A)|_{D(I)} =D(IA) \rightarrow D(I)$ is (Zariski-)locally an ${\mathbb A}^1$-scheme (i.\,e., locally the rings are polynomial rings in one variable).

\item
$\Gamma(D(IA), \O_A)=A$.
\end{enumerate}
\end{definition}

\begin{remark}
The first condition means that the natural restriction homomorphism $R=\Gamma(X, \O_X) \rightarrow \Gamma(D(I), \O_X)$ is an isomorphism ($X=\Spec R$). The second condition means that there exist ideal generators $(f_1, \ldots , f_n)=I$ such that $A_{f_i} \cong R_{f_i} [U]$ (we do not impose any condition on  the transition maps). If $R$ and $A$ are domains, then the third condition means that $A_{f_1} \cap \ldots \cap A_{f_n}= A$.
\end{remark}

\begin{proposition}
Let $R$ be a commutative ring, $I \subseteq R$ an ideal fulfilling $\Gamma(D(I),\O_X)=R$ and let $c \in H^1(U,\O_X)$ with corresponding torsor $T=T(c)$ over $U=D(I)$. Then $\Gamma(T, \O_{T})$ is an ${\mathbb A}^1$-patch in the sense of Definition \ref{a1patchgeneral}.
\end{proposition}
\begin{proof}
First of all, $T$ is quasiaffine, since it is affine over the quasiaffine scheme $D(I)$ (\cite[Corollaire 5.1.8]{EGAII} or use the explicit description of Lemma \ref{cechforcinglemma} below). For a quasiaffine scheme $T$ the natural morphism $T \rightarrow \Spec \Gamma(T,\O_T)$ is an open immersion. Hence it induces an isomorphism between $T$ and an open subset $D(J) \subseteq \Spec \Gamma(T,\O_T) $, and $\Gamma(D(J),\O_T) = \Gamma(T,\O_T)$. Furthermore, by the universal property of the ring of global sections, $\Gamma(T,\O_T)$ is an $R$-algebra. So we have a commutative diagram
\[  \begin{matrix} T & \cong  & D(J) & \subseteq & \Spec \Gamma(T, {\mathcal O}_T)  \\ \! \!\!\! \pi \!\!   \downarrow &  & &  & \downarrow \\ U & & \subseteq & & \Spec R  \, .  \end{matrix} \]
For $D(f) \subseteq U$ the open subset $T_f = \pi^{-1}(D(f)) \subseteq T$ is an affine scheme and therefore we have an isomorphism \[T_f \stackrel{\cong}{ \longrightarrow} \Spec \Gamma(T_f, {\mathcal O}_T)=\Spec \Gamma(T, {\mathcal O}_T)_f =D(f) \, .\] Hence $D(J)=D(I \, \Gamma (T, {\mathcal O}_T))$ and all properties of Definition \ref{a1patchgeneral} are fulfilled.
\end{proof}

Example \ref{5dimquadric} and the Examples in Section \ref{duboulozsection} provide natural examples of ${\mathbb A}^1$-patches in the sense of Definition \ref{a1patchgeneral}, but not in the sense of Definition \ref{a1patch}.

\begin{proposition}
\label{patchgeneralrealize}
Let $R$ be a commutative ring and let $A$ be an ${\mathbb A}^1$-patch in the sense of Definition \ref{a1patchgeneral} (with respect to some ideal $I$). Then there exists a line bundle $L$ over $U=D(I)$ and an $L$-torsor $T \rightarrow U$ such that $A=\Gamma(T, {\mathcal O}_T)$.
\end{proposition}
\begin{proof}
Let $U= \bigcup_{j \in J} D(f_j)$ be an open cover such that $A_{f_j} \cong R_{f_j}[Y_j]$. We show that $\Spec A|_U$ is an ${\mathbb A}^1$-torsor over $U$. The transition mappings are given by $Y_j = \frac{a_{ij}}{f_i^mf_j^m}Y_i + \frac{b_{ij}}{f_i^mf_j^m}$, where the $\frac{a_{ij}}{f_i^mf_j^m} $ are units in $R_{f_if_j}$. Because these data stem from an algebra, they fulfill the cocycle condition and therefore we get a cohomomology class $c \in H^1(U,{\mathcal O}^\times_X) \cong \Pic U$. Let $\shL$ be the corresponding invertible sheaf and let $L \rightarrow U$ be the corresponding line bundle which is given as $L=\Spec \bigoplus_{k \in \N} \shL^k$. Again, $L$ is given locally as the affine line over the base, and the transition mappings are by construction given by $Z_j =\frac{a_{ij}}{f_i^mf_j^m}Z_i $. We define locally 
\[ R_{f_j} [Y_j] \longmapsto R_{f_j}[Y_j,Z_j],\, Y_j \longmapsto Y_j +  Z_j  \, . \]
This action of the affine line is compatible with the transition maps for $ \Spec A$ and for $L$, hence we get an action
\[  L \times \Spec A|_U \longrightarrow \Spec A|_U \, \]
which exhibits $\Spec A|_U$ as an $L$-torsor.
\end{proof}

\begin{corollary}
Let $R$ be a commutative ring and let $A$ be an ${\mathbb A}^1$-patch in the sense of Definition \ref{a1patchgeneral} (with respect to some ideal $I$). Suppose that $R$ is a direct summand of $A$. Then $A=\Sym \shL$.
\end{corollary}
\begin{proof}
Let $T \rightarrow U=D(I)\subseteq X =\Spec R$ be the $L$-torsor defined by $A$ as explained in Proposition \ref{patchgeneralrealize} and let $c \in H^1(U, \shL)$ be the corresponding cohomology class. We have to show that $c=0$, for then $T \cong L $ and hence \[ A=\Gamma (L, {\mathcal O}_L) = \bigoplus_{n \in \N} \Gamma (U, {\shL}^{\otimes n } ) \, \]  
Assume that $c \neq 0$. Now, the pull-back of a torsor to itself is always the trivial torsor, hence $\pi^{*}(c)=0$ in $H^1(T,{\mathcal O}_T)$. We use the identification $H^1(U,\shL) \cong H^2_{X \setminus U} (M)$, where $M$ is an $R$-module which induces the invertible sheaf $\shL$ on $U$. Hence the map $R \rightarrow A$ annihilates non-trivial local cohomology. But since $A=R \oplus V$, this is not possible.
\end{proof}

\section{Background on tight closure, solid closure and affineness}

\label{backgroundsection}

We are interested in the question when a torsor over a given base scheme is an affine scheme. This is a difficult question in general. Surprisingly, the theory of \emph{tight closure}, a closure operation in positive characteristic introduced by Hochster and Huneke \cite{hochsterhunekebriancon}, \cite{hunekeapplication}, and in particular the interpretation of tight closure as solid closure provides a new way to look at this question. Let $R${} be a noetherian domain of positive characteristic $p$, let
\mathdisp {F:R \longrightarrow R,\, f \longmapsto f^p, \,} {}
be the \betonung{Frobenius homomorphism}, and let
\mathdisp {F^e:R \longrightarrow R,\, f \longmapsto f^q, \, q=p^e \, ,} {}
be its $e${}th iteration. Let $I =(f_1, \ldots ,f_n)$ be an ideal and set
\mathdisp {I^{[q]}= (f_1^q, \ldots , f_n^q) = \text{ extended ideal of } I \text{ under } F^e } {.}

\begin{definition}
The
\betonung{tight closure} of $I${} is the ideal
\mathdisp {I^*:=\{f \in R:\, \text{ there exists } z \neq 0 \text{ such that } zf^q \in I^{[q]} \text{ for all } q=p^e\} \, .} {}
\end{definition}

The relation between tight closure and forcing algebras is given in the following theorem due to Hochster (combine \cite[Definition 5.1]{hochstersolid},  \cite[Theorem 8.5(i)]{hochstersolid}, \cite[Theorem 8.6]{hochstersolid} and \cite[Corollary 2.4]{hochstersolid}).

\begin{theorem}
\label{solidlocalcriterion}
Let $R${} be a domain which is essentially of finite type over an excellent local ring of positive characteristic. Let
\mathl{f_1, \ldots, f_n}{} generate an ideal $I$ and let $f${} be another element in $R$. Then
\mathl{f \in I^*}{} if and only if
\mathdisp {H^{{\mathrm ht} ({\mathfrak m}')}_{{\mathfrak m}' } (B') \neq 0} {}
for the maximal ideals ${\mathfrak m}'$ in the complete domains $R'$ of $R$ (that is, $R'$ is the reduction modulo a minimal prime ideal of the completion of a localization $R_{\fom}$ of a maximal ideal $\fom$ of $R$), where
\mathl{B'=R'[T_1, \ldots , T_n]/(f_1T_1 + \ldots +f_nT_n +f)}{} denotes the forcing algebra of these elements. 
\end{theorem}

For $\fom$-primary ideals in a normal local domain $(R, \fom)$ this criterion becomes easier.

\begin{corollary}
\label{solidcohocrit}
Let $R${} be a normal excellent local domain with maximal ideal $\mathfrak m${} over a field of positive characteristic. Let
\mathl{f_1, \ldots, f_n}{} generate an $\mathfrak m${-}primary ideal $I${} and let $f${} be another element in $R${.} Then
\mathl{f \in I^*}{} if and only if
\mathdisp {H^{\dim (R)}_{\mathfrak m} (B) \neq 0} {,}
where
\mathl{B=R[T_1, \ldots , T_n]/(f_1T_1 + \ldots +f_nT_n +f)}{} denotes the forcing algebra of these elements. 
\end{corollary}

The definition of \emph{solid closure} can now be given by using the local criterion of Theorem \ref{solidlocalcriterion}. This gives a closure operation for commutative rings in all characteristics, denoted by $I^\star$, and Theorem \ref{solidlocalcriterion} shows that in positive characteristic under mild conditions solid closure equals tight closure.

Let us reformulate the criterion for solid closure in the situation of Corollary \ref{solidcohocrit} and relate it to the affineness of torsors. If the dimension $d={\rm dim} (R)$ is at least two, then
\mathdisp {H^d_{\mathfrak m} (R) \longrightarrow H^d_{\mathfrak m} (B) \cong H^d_{\mathfrak m B} (B) \cong H^{d-1}(D({\mathfrak m B), \mathcal O_B})} {.}
This means that we have to check the cohomological properties of the complement of the (exceptional) fiber over the closed point. In the dimension two case this yields the following.

\begin{corollary}
\label{solidtwodim}
Let $(R, \fom)$ denote a two-dimensional normal noetherian local domain, let $I=(f_1, \ldots ,f_n)$ denote an $\fom$-primary ideal and let $f \in R$ be an element. Then the following hold.

The open subset $D(\fom B) \subseteq \Spec B$ (where $B$ denotes the forcing algebra for these data) is an affine scheme if and only if $f \not\in I^\star$.

If $I=(f_1,f_2)$, then the torsor over $D(\fom)$ given by the ${\check{C}}$ech cohomology class $f/f_1f_2$ is an affine scheme if and only if $ f \not\in (f_1,f_2)^\star$.
\end{corollary}
\begin{proof}
Since the dimension is two, according to the cohomological criterion for solid closure given in Corollary \ref{solidcohocrit} we have to look whether the first sheaf cohomology of the structure sheaf vanishes. This is true if and only if the open subset \mathl{U=D(\fom B)}{} is an affine scheme by the cohomological criterion of Serre and since $U$ is quasiaffine.

The second statement follows from the first and Proposition \ref{restrictedforcing} (and the begin of Section \ref{fibrationsection}). 
\end{proof}

\begin{example}
Let $(R, \fom)$ be a two-dimensional local ring and let $T$ be the trivial $\Syz(f_1, \ldots ,f_n)$-torsor (for some $\fom$-primary ideal $(f_1, \ldots, f_n)$) over $U= \Spec R \setminus \{\fom \}$ corresponding to the cohomology class $0 \in H^1(U,\Syz(f_1, \ldots ,f_n))$. Then the torsor is just the vector bundle itself. Its total space is not an affine scheme, since there exists the zero-section of the bundle which gives a closed subscheme isomorphic to the punctured base spectrum $\Spec R \setminus \{ \fom \}$, which is not affine.
\end{example}

Results of tight closure theory give only results on the affineness of torsors in positive characteristic. In arbitrary characteristics, solid closure is strictly speaking only a reformulation of the affineness problem, but still this approach gives a good idea what to expect. 

Since a regular ring in positive characteristic is $F$-\emph{regular}, meaning that $I=I^*$ for all ideals, it follows that the induced torsor $T \rightarrow U=D(I)$ has for every prime ideal $\fop$ of height two the property that $T_{\fop} \rightarrow U_\fop=U \cap \Spec R_\fop$ is affine if and only if it is not trivial. In paricular, for a local regular ring in positive characteristic of dimension two, a torsor over the punctured spectrum is either the trivial torsor or an affine scheme. Over arbitrary open subschemes of a regular ring it is not true that every non-trivial torsor is an affine scheme, as the following two easy examples show.

\begin{example}
We consider a two-dimensional regular nonlocal ring and remove two closed points from its spectrum, e.\,g.\  $U= {\mathbb A}^2  \setminus \{ (x,y), (x,y-1)\}$, and consider the cohomology class $\frac{1}{xy}$. The torsor corresponding to this class is not an affine scheme, because the affine base-change to the localization $K[x,y]_{(x,y-1)}$ yields the trivial torsor. On the other hand, the localization to $K[x,y]_{(x,y)}$ shows that this torsor is not trivial.
\end{example}

\begin{example}
Consider a regular local ring $R$ of dimension three and let $\fom =(x,y,z)$ be its maximal ideal, e.\,g.\ $R=K[x,y,z]_{(x,y,z)}$ or $R=K[[x,y,z]]$. Consider the cohomology class
\[ c= \frac{z}{xy}  \in H^1 (D(x,y), {\mathcal O_R})  \]
and let $T(c) \rightarrow D(x,y)$ be the corresponding torsor. On one hand, the pull-back of this cohomology class to the two-dimensional localization $R_{(x,y)}$ is nonzero, as it is given by a unit divided by $xy$, hence $c$ itself is nonzero. On the other hand, the pull-back of this cohomology class to the closed subscheme
\[D(x,y) \cap V(z) \subseteq D(x,y) \]
is zero, therefore the restriction of $T(c)$ to the closed subscheme given by $z=0$ is not affine and so $T(c)$ can not be affine. In positive characteristic this means that $z \not\in (x,y)^*$, though $T(c)$ is not affine.

The forcing algebra for these data is
\[ K[x,y,z][u,v]/(ux+yv+z) \cong K[x,y,u,v] \, .\]
From this perspective it can also be immediately seen that the torsor is not affine, as it is isomorphic to  $D(x,y) \subseteq {\mathbb A}^4$. Note that the forcing algebra in this example is just a polynomial ring in four variables. This example occurs also in \cite[Example 6.2]{bonnetsurjectivity}.
\end{example}

For regular rings of dimension two not containing a field of positive characteristic, the theory of tight closure does not immediately give an answer whether the torsors are affine or not. However, the notion of solid closure gives a clear hint what to expect. In particular, we have the theorem due to Hochster that regular rings of dimension two are (weakly) S-regular \cite[Theorem 7.20]{hochstersolid}, meaning that every ideal equals its solid closure (this is an application of the monomial conjecture which holds in dimension two in all characteristics).

\section{Affineness and global generation of torsors}

\label{affinesection}

We deal now with the case of a two-dimensional local base ring (of any characteristic) and the ${\mathbb A}^1$-torsors on the punctured spectrum. In particular, we investigate the question when every non-trivial torsor is an affine scheme. We start with the regular case.
\begin{lemma}
\label{regularaffinea1}
Let $R$ be a two-dimensional regular local ring. Then every non-trivial ${\mathbb A}^1$-torsor on $U=\Spec R \setminus \{\fom\} $ is an affine scheme.
\end{lemma}
\begin{proof}
An ${\mathbb A}^1$-torsor $T$ over $U$ is given by a cohomology class $c \in H^1(U, \O)$. This class is non-zero, because the bundle is non-trivial. Let $c= \frac{a}{x^n y^m}$ be a ${\rm \check{C}}$ech representation with regular parameters $x,y$ and $a \not \in (x^n,y^m)$ (else the class would be $0$). Since $R$ is S-regular by \cite[Theorem 7.20]{hochstersolid}, it follows that $a \not\in (x^n,y^m)^\star$ and from Corollary \ref{solidtwodim} we deduce that $T$ is an affine scheme.
\end{proof}

\begin{remark}
For the affine plane this statement is also proved in \cite[Proposition 1.2]{duboulozfinstonexotic}). In the complex situation, M. Abe has shown that any non-trivial algebraic ${\mathbb C}$-torsor over the punctured plane is a Stein manifold, \cite[Theorem 3.1]{abetotalspace2}. These results follow from Lemma \ref{regularaffinea1}.
\end{remark}

\begin{theorem}
\label{regularaffine}
Let $R$ be a two-dimensional regular local ring. Then every non-trivial affine-linear bundle on $U=\Spec R \setminus \{\fom\} $ is an affine scheme.
\end{theorem}
\begin{proof}
An affine-linear bundle over $U$ is given by a vector bundle $\shS$ on $U$ and a cohomology class $c \in H^1(U, \shS)$. Since we suppose that the bundle is non-trivial, this class is nonzero. By Hilbert's syzygy theorem, $\shS$ is a free $R$-module, say $\shS=R^n$. Let $c=(c_1 , \ldots , c_n)$. The torsor corresponding to this class is $T=T_1 \times_U T_2 \times_U \cdots \times_U T_n $, where $T_i$ is the torsor corresponding to $c_i$. The morphisms $T_i \rightarrow U$ are affine morphisms, hence all projections from $T $ to a product with less components are also affine morphisms. So if at least one $T_i$ is affine, $T$ must be affine. But since $c_i \neq 0$ for at least one $i$, we get the affineness by Lemma \ref{regularaffinea1}.
\end{proof}

Recall that a quotient singularity is by definition the spectrum of the invariant ring of a regular ring on which a finite group is acting (with group order not divisible by the characteristic).

\begin{theorem}
\label{quotientaffine}
Let $R$ be a two-dimensional local quotient singularity over a field $K$. Then every non-trivial affine-linear bundle on $U=\Spec R \setminus \{\fom\} $ is an affine scheme.
\end{theorem}
\begin{proof}
Let $R \subseteq S$ be a finite extension such that $R=S^G$ with $(S, \fon)$ regular local and $G$ a finite group of order prime to the characteristic. Then the trace map divided by the group order shows that $R$ is a direct summand of $S$. An affine-linear bundle over $U$ is given by a vector bundle $\shS$ on $U$ and a cohomology class $c \in H^1(U, \shS)$. This class is non-zero, because the bundle is non-trivial. The pull-back of this class to $V=\Spec S \setminus \{\fon \}$ is still non-zero, because $R$ is a direct summand of $S$. The product $T'=V \times_U T$ is the torsor given by $\varphi^*(c)$ and this is finite over $T$. By Theorem \ref{regularaffine}, it is an affine scheme. Therefore by the theorem of Chevalley \cite[Theorem 6.7.1]{EGAII}, $T$ itself is an affine scheme.
\end{proof}

We do not know whether this statement also holds for a two-dimensional rational singularity in characteristic zero (see below).

We now answer the last question of \cite{duttaguptaonoda}. It asks whether an ${\mathbb A}^1$-patch $A$ (in the sense of Definition \ref{a1patch}) over a factorial two-dimensional domain $R$ is finitely generated over $R$, in particular for the icosahedral singularity. We show that for this special case the answer is yes, but not for factorial domains in general.

\begin{example}
Let $R=K[X,Y,Z]/(X^2+Y^3+Z^5)$ and suppose that the characteristic is at least $7$. Then $R$ is a quotient singularity, so in particular $R$ is a direct summand of $K[U,V]$. Hence the affineness of any non-trivial ${\mathbb A}^1$-torsor (and flatness and finite generation) follows from Theorem \ref{quotientaffine} and Proposition \ref{patchproperties}(iii).

For example, the torsor induced by the forcing equation $YU+ZV+X$ (and the cohomology class $\frac{X}{YZ}$), i.\,e.\ the open subset
\[ D(Y,Z) \subset \Spec R[U,V]/( YU+ZV+X )\]
is an affine scheme, just because $X \not\in (Y,Z)$. This may fail in small characteristics. For example, in characteristic two, we have the Frobenius inclusion $X^2 \in (Y^2,Z^2)$, hence the torsor trivializes after applying the Frobenius and cannot be affine.
\end{example}

\begin{remark}
The techniques developed in \cite[Section 8]{brennerparasolid} allow us to deduce the corresponding statements also in mixed characteristic (excluding finitely many residue characteristics). E.g. if $S$ is a three-dimensional regular local ring of mixed characteristic $p$ with maximal ideal $(p,y,z)$, then $R=S/(p^2+y^3+z^5)$ is pararegular \cite[Definition 4.2]{brennerparasolid} for $p \gg 0$ and the torsor given by the cohomology class $\frac{p}{yz}$ is affine as well. In particular, the subset
\[D(Y,Z) \subseteq  \Spec {\mathbb Z} [Y,Z]_{(p,Y,Z)}[U,V]/(YU+ZV+p)\]
is affine for $p\gg 0$.
\end{remark}

\begin{example}
Let $R=K[X,Y,Z]/(X^2+Y^3+Z^7)$, which is a factorial domain by \cite{samuelfactorial}. We consider again $X$, which does not belong to $(Y,Z)$. $R$ can be graded by assigning $\deg (X)=21$, $\deg (Y) = 14$ and $\deg (Z)=6$. Hence the corresponding ${\rm \check{C}}$ech cohomology class $\frac{X}{YZ}$ has degree $21-14-6=1$, which is positive (in particular the $a$-invariant is not negative and by a Theorem of Flenner and Watanabe they do not have a rational singularity  \cite{flennerrationalquasihomogen}, \cite{watanaberational}).

Suppose first that we are in positive characteristic $p$, and assume that the characteristic is large enough such that $R$ is normal. Then for a sufficiently high power $q=p^{e}$, the Frobenius pull-back $\frac{X^q}{Y^qZ^q}$ of the cohomology class has arbitrarily large degree and must therefore vanish. Hence $X^q \in (Y^q,Z^q)$ and so we have a Frobenius-trivialization of the torsor which can therefore not be affine (and its ring of global sections is not flat; finite generation is not known).

Suppose now that we are in characteristic zero. Since affineness is not affected by flat base field extensions we may assume that $K= {\mathbb Q}$. We look at the family
\[  D(Y,Z)  \subset \Spec  {\mathbb Z}[X,Y,Z][U,V]/(X^2+Y^3+Z^7, YU+ZV+X ) \longrightarrow \Spec {\mathbb Z} \, .  \]
If the fiber over the generic point was affine, then almost all fibers would be affine, but this is not true by the consideration in positive characteristic. Hence the torsor $D(Y,Z)$ is not affine in characteristic zero. By \cite[Corollary 1.6]{brennertightproj}, the extended ideal $(Y,Z)$ in the forcing algebra has superheight one. Its ring of global sections is therefore not finitely generated by \cite[Theorem 3.2]{brennersuperheight}.
\end{example}

In the following statement we use the notion of a \emph{rational singularity}. For this notion and its main properties see \cite{lipmanrational} or \cite{bingenerstorchdivisor}, for $F$-rationality and its relation to rational singularities see \cite[Section 10.3]{brunsherzog}, \cite{smithrational}, \cite{hararationalfrobenius}.

\begin{theorem}
\label{nonrationalnonaffine}
Let $(R, \fom) $ be a two-dimensional excellent normal local domain over a field $K$ which is not a rational singularity. Then there exists a non-trivial ${\mathbb A}^1$-torsor on $U=\Spec R \setminus \{\fom\} $ which is not an affine scheme. If the characteristic of $K$ is zero, then there exists an ${\mathbb A}^1$-torsor on $U $ which is not an affine scheme and such that its ring of global sections is not finitely generated over $R$.
\end{theorem}
\begin{proof}
Suppose first that the characteristic of $R$ is positive. Then $R$ is not $F$-rational by a Theorem of K. Smith \cite{smithrational}. This means that there exists a system of parameters $f, g \in R$ such that $(f,g)^* \neq (f,g)$. Let $h \in (f,g)^* \setminus (f,g)$. Then because of $(f,g)^* =(f,g)^\star$ (tight closure is solid closure), we know that the torsor given by the cohomology class $h/fg$ is not trivial, but also not affine.

Suppose now that the characteristic of $K$ is zero. Set $X= \Spec R$ and let $\tilde{X} \rightarrow X$ be a resolution of singularities and denote the exceptional fiber (the fiber over the maximal ideal $\fom$) by $E$. Since $X$ does  not have a rational singularity there exists a nonzero cohomology class $ \frac{h}{fg} = c \in H^2_{\fom} (R) \cong H^1(U, \O_X)$ which becomes $0$ in $H^2_E( \tilde{X})$. Equivalently, this means that $c \in H^1(U, \O_{\tilde{X} })=  H^1(U, \O_X)$ stems from a cohomology class $\tilde{c} \in H^1(\tilde{X}, \O_{\tilde{X } })$ by restriction.

We show that the torsor $T(c)$ on $U=D(f,g)$ is not an affine scheme. This torsor equals the open subset
\[D(f,g) \subseteq \Spec\, ( R[T_1,T_2]/(fT_1+gT_2+h) ) \, \]
as explained in Proposition \ref{restrictedforcing}. The superheight of $(f,g)A$ is $1$ by \cite[Corollary 1.6]{brennertightproj}. Hence from the non-affineness we get by \cite[Theorem 3.2]{brennersuperheight} that $\Gamma(T(c), \O_{T(c)})$ is not finitely generated over the forcing algebra nor over $R$.

So assume that $T(c)$ is an affine scheme. Then we can express all relevant data and properties in a finitely generated ${\mathbb Z}$-algebra $B \subseteq R$ including normality, that the radical of $(f,g)$ is a prime ideal of height $2$, the non-zero cohomology class, the vanishing of it in the resolution, the affineness of the torsor. Therefore we may assume that $R$ is essentially of finite type over a finitely generated $\mathbb Z$-algebra. But then the affineness of the torsor descends to almost all prime reductions producing a contradiction to the statement in positive characteristic.
\end{proof}

\begin{remark}
We do not know whether in the situation of Theorem \ref{nonrationalnonaffine} the ring of global sections of the non-trivial non-affine torsor in positive characteristic is finitely generated or not. Nor do we know whether the converse of the statement holds in characteristic zero, namely whether in characteristic zero a two-dimensional rational singularity implies that the only non-affine ${\mathbb A}^1$-torsor is the trivial one. For a quotient singularity this was proved in Theorem \ref{quotientaffine} above. In positive characteristic it follows from the equivalence of rational singularities and $F$-rational singularities proved by K. Smith \cite{smithrational} and N. Hara \cite{hararationalfrobenius}.
\end{remark}

\section{Forcing algebras and ${\rm \check{C}}$ech cohomology classes}
\label{cechsection}

In this section we describe another situation where forcing algebras arise naturally. Namely by annihilating a ${\rm \check{C}}$ech cohomology class by expressing a ${\rm \check{C}}$ech cocycle as a coboundary with the help of new indeterminates. This can be studied much more generally, but we restrict to  first cohomology of the structure sheaf in a quasiaffine scheme. For ${\rm \check{C}}$ech cohomology in general we refer to \cite[Chapter III.4]{hartshornealgebraic}.

Let $U=D(f_1, \ldots ,f_n) \subseteq \Spec \, R$. Then a cohomology class $c \in H^1(U, \O_U)$ is represented in the following way. Take
\[ \beta = ( \frac{ \classnen_{ij}}{f_i^mf_j^m} )_{i,j} \]
fulfilling the cocycle condition, namely
\[ (d \beta)_{ijk} = \frac{ \classnen_{ij}}{f_i^mf_j^m} -\frac{ \classnen_{ik}}{f_i^mf_k^m} + \frac{ \classnen_{jk}}{f_j^mf_i^m} = 0 \, \]
for all $i,j,k$, or, equivalently if $R$ is a domain,
\[ \classnen_{ij}f_k^m - \classnen_{ik} f_j^m +  \classnen_{jk} f_i^m =0 \, \]
The condition whether this cohomology class is zero, meaning that its representative is in the image of the coboundary map, can be expressed with a matrix of $n$ columns and $\binom{n}{2}$ rows in the following way.
\[   \begin{pmatrix} \frac{1}{f_1^m}  & -\frac{1}{f_2^m} & 0 & \ldots & 0  \\  \frac{1}{f_1^m}  & 0 & -\frac{1}{f_3^m} &  \ldots & 0 \\ \vdots & \vdots & \vdots &  \ddots &\vdots \\ 0 & \ldots &0&  \frac{1}{f_{n-1}^m}  & - \frac{1}{f_n^m}   \end{pmatrix}  \begin{pmatrix} t_1 \\ t_2 \\ \vdots \\ t_n   \end{pmatrix} =   \begin{pmatrix} \frac{\classnen_{12}}{f_1^mf_2^m}  \\  \frac{\classnen_{13}}{f_1^mf_3^m}  \\ \vdots \\ \frac{\classnen_{n-1 n }}{f_{n-1}^m f_n^m}   \end{pmatrix} \]
(it might be necessary to increase $m$, but we work with such systems for all possible exponents anyway).
We multiply this in each $(i,j)$-row with $f_i^mf_j^m$ and get
\[  \begin{pmatrix} {f_2^m}  & -{f_1^m} & 0 & \ldots & 0  \\  {f_3^m}  & 0 & -{f_1^m} &  \ldots & 0 \\ \vdots & \vdots & \vdots &  \ddots &\vdots \\ 0 & \ldots &0&  f_{n}^m  & -f_{n-1}^m \end{pmatrix}  \begin{pmatrix} t_1 \\ t_2 \\ \vdots \\ t_n   \end{pmatrix} =   \begin{pmatrix} \classnen_{12}  \\  \classnen_{13} \\ \vdots \\ \classnen_{n-1 n }  \end{pmatrix}   \, .\]
This is a system of forcing equations, which we write briefly as $\repmatrix t= s$. The corresponding forcing algebras depend of course on the exponent $m$, but the induced scheme over $U$ is an ${\mathbb A}^1$-torsor which depends only on the cohomology class, as the following Lemma shows.

\begin{lemma}
\label{cechforcinglemma}
Let $R$ be a commutative ring, $f_1 , \ldots, f_n \in R$, $U=D(f_1, \ldots, f_n)$ and let $c \in H^1(U, {\mathcal O_X})$ be represented by the ${\check{C}}$ech cocycle $c= (\frac{\classnen_{ij}}{f_i^mf_j^m} )_{i,j}$. Then the forcing algebra
\[B=R[T_1 , \ldots ,T_n]/(\repmatrix T- s)\]
has the property that the restriction
\[ (\Spec B)|_{U}\cong T(c) \]
is isomorphic to the torsor given by $c$.
\end{lemma}
\begin{proof}
We consider the localizations $B_{f_i}$. With the help of the $(i,j)$-row ($j \neq i$) we can eliminate every $T_j$ ($j \neq i$) and therefore we get
\[ B_{f_i} \cong R_{f_i}[T_i] \cong R_{f_i}[W_i] \, , \]
where we set $W_i= \frac{T_i}{f_i^m}$. For the $W_i$ we have the transformation rule
\[W_i =W_j + \frac{\classnen_{ij}}{f_i^mf_j^m} \, ,\]
so this is a realization of the cohomology class as a torsor.
\end{proof}

\begin{remark}
\label{cechthree}
For three ideal generators $f_1,f_2,f_3$ the representation is
\[ ( \frac{b_{3}}{f_1^mf_2^m} ,   \frac{b_{2}}{f_1^mf_3^m} , \frac{b_{1}}{f_2^mf_3^m}  )  \]
with the condition that
\[ b_3 f_3^m - b_2f^m_2+b_1f^m_1 =0  \, .\] 
This gives the system of forcing equations
\[  \begin{pmatrix}  f_2^m & - f_1^m & 0   \\ f_3^m & 0 & -f_1^m \\ 0 & f_3^m & -f_2^m   \end{pmatrix} \begin{pmatrix} T_1 \\ T_2 \\ T_3  \end{pmatrix} = \begin{pmatrix} b_3 \\ b_2 \\ b_1   \end{pmatrix}\, . \]
\end{remark}

\begin{example}
\label{5dimquadric}
Let
\[R=K[X,Y,Z,U,V,W]/(UX+VY+WZ)\]
and consider the ${\rm \check{C}}$ech cohomology class
\[ ( \frac{W}{XY} , -  \frac{V}{XZ} , \frac{U}{YZ}  )  \]
on $D(X,Y,Z)$. The ring equation shows that the cocycle condition is fulfilled. The forcing equations for this are
\[   \begin{pmatrix}  Y & - X & 0   \\ Z & 0 & -X \\ 0 & Z & -Y   \end{pmatrix} \begin{pmatrix} T_1 \\ T_2 \\ T_3  \end{pmatrix} = \begin{pmatrix} W \\ -V \\ U   \end{pmatrix} \, .\]
The induced torsor is not affine, because the restriction to the affine plane given by $X=U=V=W=0$ gives the trivial, non-affine torsor over the punctured plane. The torsor is also not trivial, as the restriction to the punctured plane given by $Z=U=V=0$ and $W=1$ shows.
\end{example}

For other examples of forcing algebras coming from a ${\rm \check{C}}$ech cohomology class  see Section \ref{duboulozsection}.

\section{Torsors over the punctured plane and the cancellation problem}
\label{cancellationsection}

For the following argument compare also the introduction of \cite{finstonmaubach}.

\begin{lemma}
Let $U$ be a separated scheme and let $c,c' \in H^1(U, {\mathcal O}_U)$ be cohomology classes with corresponding torsors $T(c)$ and $T(c')$. Suppose that $T(c)$ is an affine scheme. Then \[T(c) \times_U T(c') \cong T(c) \times_U {\mathbb A}^1_U \, . \]

If $T(c)$ and $T(c')$ are both affine, then 
\[ T(c) \times {\mathbb A}^1  \cong \,   T(c') \times_U {\mathbb A}^1_U \,   . \]
\end{lemma}
\begin{proof}
We consider the commutative diagram
\[  \begin{matrix}  & & T(c) \times_U T(c') & &   \\ & \swarrow & & \searrow &  \\ T(c) & & & & T(c') \\ & \searrow & & \swarrow & \\ & & U & &\end{matrix}  \, .\]
The arrows in the first row are the pull-backs of the arrows in the second row, and so they are the morphisms to the base scheme of the pull-backs of the cohomology classes. In particular, $T(c) \times_U T(c') \cong  T(p^*(c'))$, where $p:T(c) \rightarrow U$. If $T(c)$ is affine, then $p^*(c')=0$ is trivial and hence $T(p^*(c')) = T(c) \times_U {\mathbb A}^1_U$.
\end{proof}

If $U$ is a separated scheme over a field $K$, then we will also write $  T(c) \times {\mathbb A}^1$ for the product over the base scheme $\Spec K$. With this lemma one can get interesting candidates for counterexamples for the general cancellation problem, namely schemes $T$ and $T'$ with $T \times {\mathbb A}^1 \cong T' \times {\mathbb A}^1$, but $T \not\cong T'$. If $T(c)$ is affine, then there is a $U$-morphism $T(c) \rightarrow T(c')$ (since the pull-back of $T(c')$ has a section); so if both are affine, then there exist $U$-morphisms in both directions, but they are not necessarily invers to each other.

\begin{example}
\label{torsoraffineplane}
We consider the family ($m,n \geq 1$)
\[B_{m,n} =K[x,y][t_1,t_2](x^mt_1 +y^nt_2-1) \, ,\]
which are the rings for the affine torsors corresponding to the cohomology classes $ \frac{1}{x^ny^m}$. Note that these forcing algebras have empty fiber over the maximal ideal, so its spectrum is already the torsor which we denote by $X_{m,n}$.  It is clear that they are not isomorphic as torsors over the punctured plane. Topologically, they are all homeomorphic to ${ \mathbb C}^3 \setminus \{\mathrm{a \, \, line} \}$. Because of their affineness we have $X_{m,n} \times {\mathbb A}^1 \cong X_{p,q}  \times {\mathbb A}^1 $ for all possible values of the indices. It was shown in 
\cite[Theorem 1]{duboulozfinstonmehta} (see also \cite[Example 2.4]{duboulozfinstonexotic}) that $X_{m,n}  \cong X_{p,q}  $ in case $m+n = p+q$. It was shown in \cite[Corollary 2.6]{duboulozfinstonexotic} that $X_{m,n}$ for $m+n \geq 3$ is not isomorphic to $Sl_2 = X_{1,1}$.
\end{example}

\section{A class of examples}
\label{duboulozsection}

The following example of a scheme (and the torsors over it) was brought to my attention by A. Dubouloz at the Bangalore conference. This is also directly related to an example of J. Winkelmann \cite[Section 2]{winkelmannfreeholomorphic}. Let
\[R=K[x,y,u,v,z]/(xv+yu+z(z-1)) \, \]
over an arbitrary field $K$ of characteristic $p \geq 0$. The partial derivatives are
\[ (v,u,2z-1, y,x) \, ,\]
so its spectrum is a smooth four-dimensional affine variety (consider the cases $p=2$ and $p \neq 2$ separately). We consider the complement of $F=V(x,y,z)$ which is a subset of codimension two. What can we say about the ${\mathbb A}^1$-torsors over
\[U=D(x,y,z)= \Spec R \setminus F \, .\]
This question is interesting with respect for the affine cancellation problem, i.\,e.\_ the question whether ${\mathbb A}^{n+1} \cong X \times {\mathbb A}^1$ implies ${\mathbb A}^{n}$.

\begin{example}
\label{examplequadrictorsor}
We consider on $U$ the first ${\rm \check{C}}$ech cohomology class (compare Remark \ref{cechthree})
\[ c= (\frac{z-1}{xy} , - \frac{u}{xz}  , \frac{v}{yz} ) \, .\]
The equation of the variety shows immediately that this is a cocycle. The system of forcing equations of this class is
\[  \begin{pmatrix}  y & - x & 0  \\ z & 0 & - x \\  0 & z & -y  \end{pmatrix} \begin{pmatrix} t_1 \\ t_2 \\  t_3  \end{pmatrix} = \begin{pmatrix} z-1 \\ -u \\ v   \end{pmatrix} \, .\]
The first row shows that  $z-1$ belongs to the ideal generated by $x$ and $y$, hence $x,y,z$ generate the unit ideal in the forcing algebra. So the fibers of the spectrum of this forcing algebra over $V(x,y,z)$ are empty and the torsor is immediately an affine scheme.

We can omit the basic equation $xv+yu+z(z-1)=0$ because we can reconstruct it from the three equations which appear in the above system by
\begin{eqnarray*}
 xv+yu &=& x(z t_2 -y t_3) - y(zt_1-xt_3) \\
&=& xzt_2-yz t_1 \\
&=& -z(yt_1-xt_2) \\
&=&-z(z-1) \, .
\end{eqnarray*}

Moreover, by the last row we can eliminate $v$, by the second row we can eliminate $u$ and by the first row we can eliminate $z$. Hence this forcing algebra is just the polynomial ring in the five variables $x,y,t_1,t_2,t_3$ and  the torsor $T(c)$ is isomorphic to ${\mathbb A}^5$. The mapping to $\Spec R$ sends altogether $(x,y,t_1,t_2,t_3 ) $ to
\[ (x,y,z,u,v) =(x,y,yt_1-xt_2+1, - (yt_1- xt_2+1)t_1 +xt_3,t_2 (yt_1- xt_2+1) -yt_1 ) \, .\]
\end{example}

The following example shows that not all non-trivial torsors on $U$ are affine. 

\begin{example}
\label{examplequadrictorsornotaffine}
We consider the first ${\rm \check{C}}$ech cohomology class $d=uc$, where $c$ is the class of Example \ref{examplequadrictorsor}, i.\,e.\
\[d= (\frac{u(z-1)}{xy} , - \frac{u^2}{xz}  , \frac{uv}{yz} ) \, .\]
The equation $u=0$ defines a closed subscheme 
\[V=V(u) \cong \Spec K[x,y,z,v]/( xv+z(z-1)) \subseteq \Spec R\] and $U'= U \cap V(u) \subseteq U$ is a closed subscheme of $U$ and an open subscheme of $V$. The prime ideal $(x,y,z) \subset K[x,y,z,v]/( xv+z(z-1))$ has height two, therefore $U'$ is not an affine scheme. The restriction of the cohomology class $d$ to $U'$ is $0$, since $u=0$ on $U'$, therefore $T(d)|_{U'} \cong T(d|_{U'})$ is the trivial torsor over a non-affine scheme, hence non-affine. Because $T(d)|_{U'} \subseteq T(d)$ is a closed subscheme, also $T(d)$ is not affine. We claim that $T(d)$ itself is not trivial. For this we consider the restriction of $T(d)$ to the open subset $D(x,y) \subseteq U$, which is given by the cohomology class $\frac{u(z-1)}{xy}$. We further restrict this cohomology class to the open subset $D(x,y) \subseteq \Spec S$, where $S=R_{(x,y,z,u-1,v)}$ is regular and local. There, $u(z-1)$ is a unit and hence the class $\frac{u(z-1)}{xy}$ is not trivial.
\end{example}

\begin{proposition}
Let $R=K[x,y,u,v,z]/(xv+yu+z(z-1))$ and let $c \in H^1(U, {\mathcal O}_U)$ ($U=D(x,y,z)$) be a first cohomology class with corresponding torsor $T(c)$ on $U$. Then the following hold.

(i) There exists an ${\mathbb A}^1$-torsor on $T(c)$ which is isomorphic to $ {\mathbb A}^6$. 

(ii) If $T(c)$ is affine, then $T(c) \times {\mathbb A}^1 \cong {\mathbb A}^6$. 

(iii) If $K =\CC$, then the complex space $T(c)({\mathbb C}) \rightarrow U({\mathbb C})$ is homeomorphic and $C^\infty$-diffeomorphic to ${\mathbb C}^5$.
\end{proposition}
\begin{proof}
(i). Let $c'$ be the cohomology class described in Example \ref{examplequadrictorsor} whose torsor $T(c')$ is isomorphic to ${\mathbb A}^5$. We have the commutative diagram
\[  \begin{matrix}  & & T(c) \times_U T(c')\cong {\mathbb A}^6  & &  \\ & \swarrow & & \searrow &  \\ T(c) & & & & T(c') \cong {\mathbb A}^5  \\ & \searrow & & \swarrow & \\ & & U & &\end{matrix}  \, ,\]
where the isomorphism of the fiber product with ${\mathbb A}^6 $ follows from the fact that the fiber product over the affine scheme $T(c')\cong {\mathbb A}^5$ is a trivial torsor. The projection of the fiber product to $T(c)$ is the pull-back of $T(c') \rightarrow U$ and gives an ${\mathbb A}^1$-torsor over $T(c)$.

(ii) follows from (i), since torsors over an affine scheme are trivial.

(iii). By Example \ref{examplequadrictorsor}, there exists a class $c'$ such that $T(c')$ is a fivedimensional affine space. Therefore $T(c')({\mathbb C})$ is a fivedimensional complex space and a tendimensional real space. Because of the existence of $C^\infty$-partition of unity, all torsors over a real $C^\infty$-manifold are trivial, hence the topological and the diffeomorphic shape of $T(c)$ do not depend on $c$, so they are all diffeomorphic to ${\mathbb C}^5$.
\end{proof}

By taking the trivial torsor $U \times {\mathbb A}^1$ over $U$ and on this the affine torsor coming from the pull-back of the torsor of Example \ref{examplequadrictorsor}, we see that there exists a (simply transitive) action of $({\mathbb A}^2,+)$ on the sixdimensional affine space such that the quotient is $U$ (and not an affine scheme). This group action was directly described in \cite[Section 2]{winkelmannfreeholomorphic}.

In the following example we describe a sequence of affine torsors $T(c_k)$ over $U$. We do not know if their total spaces are isomorphic to ${\mathbb A}^5$ or not, though $T(c_k) \times {\mathbb A} = {\mathbb A}^6$ is known.

\begin{example}
We consider the $k$th power of the equation $  xv +yu = z(1-z) $, i.\,e.\
\[  z^k(1-z)^k = (xv +yu)^k = x^kv^k +y( \sum_{i= 0}^{k-1}   \binom{k}{i} x^{i}v^{i} y^{k-i-1} u^{k-i} ) \, . \]
This gives rise to the cohomology class
\[ c_k = ( \frac{- (1-z)^k}{x^ky}, \frac{ - \sum_{i= 0}^{k-1}   \binom{k}{i} x^{i}v^{i} y^{k-i-1} u^{k-i} }{x^kz^k}  , \frac{ v^k}{yz^k}  ) \]
and to the system of forcing equations
\[ \begin{pmatrix}  y & -x^k & 0 \\ z^k & 0 & -x^k \\ 0 & z^k & -y \end{pmatrix} \begin{pmatrix} t_1 \\ t_2 \\ t_3 \end{pmatrix} = \begin{pmatrix} - (1-z)^k \\ - \sum_{i= 0}^{k-1}  \binom{k}{i} x^{i}v^{i} y^{k-i-1} u^{k-i}   \\ v^k \end{pmatrix}  \, \]
(the equation $xv+yu= z(1-z)$ is still around). Over the closed subset $V(x,y,z)$ these torsors have empty fibers, since then the first row has no solution. Therefore the spectra of these forcing algebras are the torsors over $D(x,y,z)$ of these cohomology classes and these torsors are all affine. We conjecture that for $k \geq 2$ these forcing algebras are not isomorphic to the polynomial ring in five variables, so these would give counterexamples to Zariskis affine space cancellation conjecture.
\end{example}

We study the forcing algebra of the previous construction for $k=2$ in more detail.

\begin{example}
We consider the algebra $B$ over $K[x,y,u,v,z]/(xv+yu+z(z-1))$ which is given by the system of forcing equations
\[ \begin{pmatrix}  y & -x^2 & 0 \\ z^2 & 0 & -x^2 \\ 0 & z^2 & -y \end{pmatrix} \begin{pmatrix} t_1 \\ t_2 \\ t_3 \end{pmatrix} = \begin{pmatrix} -(z-1)^2 \\ -  yu^2 - 2xvu   \\ v^2 \end{pmatrix}  \, .\]
With the help of the basic equation and the first equation we can eliminate $z$. By subtracting the two equations
\[ yt_1-x^2t_2 = -(z-1)^2 =-z^2+2z-1 \]
and
\[  xv+yu = z(1-z)=z-z^2\]
we get
\[ yt_1 - x^2t_2 -xv-yu = z -1 \]
 and hence
\[  z= 1  +yt_1 - x^2t_2  -xv-yu \, .  \]
Since we only have made a subtraction, we can get rid of the basic equation (or of the first equation, but  not of both, as a CoCoA-computation shows) and transform the other equations by replacing $z$. But still it seems to be quite difficult to understand this algebra, yet alone to decide whether it is a polynomial algebra or not.

The algebra $B$ can be bigraded by assigning the degrees
\[ \deg(x) = (1,0),\,  \deg(y) = (0,1),\,  \deg(z) = (0,0),\, \]
\[ \deg(v) = (-1,0),\,  \deg(u) = (0,-1),\, \]
\[ \deg(t_1) = (0,-1),\,  \deg(t_2) = (-2,0),\, \deg(t_3) = (-2,-1) \, .\]
The basic equation has degree $(0,0)$ and the forcing equations have degree $(0,0),\, (0,-1)$ and $(-2,0)$ respectively. Typical examples in the degree $(0,0)$-ring are $z,xv,yu,yt_1,x^2t_2,x^2yt_3$. How does this ring look like?

We have a closer look at the first forcing equation. We write it as
\begin{eqnarray*}
 -yt_1+x^2t_2 &=&  (y(t_1-u)-x(xt_2+v))^2 \\
& = & y^2(t_1-u)^2 + x^2(xt_2+v)^2 -2xy(t_1-u)(xt_2+v) 
\end{eqnarray*}
and
\[ x^2 (t_2 -(xt_2+v)^2) =y(t_1+y(t_1-u)^2 -2x(t_1-u)(xt_2+v)) \, . \]
Because $B$ is a factorial domain (as the polynomial ring in one variable over it is the polynomial ring in six variables over $K$), and because $x$ and $y$ are prime elements (because they are part of a system of generating variables in $B[W]$), the longer factors must be multiples of $y$ resp.\ $x^2$. We will check this explicitely. With the help of the first and the third forcing equation we compute
\begin{eqnarray*}
& & y(t_3 - t_1 t_2 +2ut_2) \\
&=& yt_3 -yt_1t_2 +2y ut_2  \\
&=&   yt_3 +yt_1t_2 -2 yt_1t_2 +2 yut_2      \\
&=& z^2 t_2 -v^2  + ( x^2t_2 -(z-1)^2) t_2  -2 yt_1t_2      +2yut_2 \\
&=& ( 2z-1 ) t_2 -v^2 + x^2t_2^2  -2 yt_1t_2   +      2yut_2 \\
&=& (1 +2yt_1 - 2x^2t_2  -2xv-2yu )t_2  -v^2 + x^2t_2^2    -2 yt_1t_2   +      2yut_2 \\
&=& t_2 -x^2t_2^2 -2xvt_2  -v^2  \, .
\end{eqnarray*}
So the above element factors as
\[  x^2y (  t_3 - t_1 t_2 +2ut_2 )  \, .\] 
\end{example}

\section{Derivations on forcing algebras}
\label{derivationsection}

We recall that a $K$-derivation $D$ on a $K$-algebra $B$ is called \emph{locally nilpotent} if for every $g \in B$ there exists a power $n$ such that $D^n(g)=0$. Locally nilpotent derivations are directly linked to group actions of the additive group. The following lemma shows that on a forcing algebra for an ideal generated by a regular sequence of two elements in characteristic zero we get a canonical locally nilpotent derivation such that its kernel is the base ring. This result is implicit in \cite{gurjarmasudamiyanishi}, at least in many examples.

\begin{lemma}
Let $R$ denote a commutative noetherian domain containing a field of characteristic zero and let $f_1,f_2, f_3 \in R$ be three elements such that the ideal $(f_1,f_2)$ has depth $2$. Then
$$D = f_2 \frac{\partial }{\partial T_1} - f_1 \frac{\partial }{\partial T_2} $$
is a locally nilpotent derivation on the forcing algebra $B=R[T_1,T_2]/(f_1T_1+f_2T_2+f_3)$ and $B^D=R$.
\end{lemma}
\begin{proof}
$D$ is a locally nilpotent derivation on $R[T_1,T_2]$ which annihilates the forcing equation, hence it is a locally nilpotent derivation on $B$.
For every $g \in R$ we have $ \frac{\partial }{\partial T_1} (g)= 0 =  \frac{\partial }{\partial T_2} (g)$, hence  $R \subseteq B^D$.

Let $g \in B^D$. The derivation $D$ is on the localization $B_{f_2} \cong R_{f_2}[ T_1] \cong R_{f_2}[W_1]$, where we set $W_1= \frac{T_1}{f_2}$, just the standard derivation with respect to the variable $W_1$. Similarly, the derivation is on $B_{f_1} \cong R_{f_1}[ T_2] \cong R_{f_1}[W_2]$ (with $W_2= \frac{T_2}{f_1}$) just the negative standard derivation with respect to $W_2$. Because we suppose characteristic $0$ we deduce $g \in R_{f_1}$ and $ g \in R_{f_2}$. Hence $g \in R$ due to the depth assumption.
\end{proof}

\section{Singular points on torsors}
\label{singularsection}

For a forcing algebra $B$ over a $K$-algebra $R=K[x_1 , \ldots , x_m]/(g_1, \ldots , g_k)$ of finite type one can determine the singular locus with the help of the Jacobian criterion. If $B$ is given as
\[B=K[x_1 , \ldots , x_m, t_1 , \ldots , t_n ]/( g_1, \ldots , g_k, \sum_{i=1}^n  f_i t_i - f  ) \, , \]
then the Jacobian matrix is (writing $h$ for the forcing equation)
\[
\begin{pmatrix}
 \frac{\partial g_1  }{\partial x_1}  & \ldots & \frac{\partial g_1 }{\partial x_m}  & 0 & \ldots & 0 \cr
\vdots & \ddots & \vdots & \vdots & \ddots& \vdots  \cr
\frac{\partial g_k }{\partial x_1}  & \ldots &\frac{\partial g_k }{\partial x_m}  & 0 & \ldots & 0 \cr
\frac{\partial h}{\partial x_1}  & \ldots  & \frac{\partial h}{\partial x_m}  & f_1 & \ldots & f_n 
\end{pmatrix} \, .
\]
Of course, $\frac{\partial h}{\partial x_j} = \sum_{i=1}^n  t_i  \frac{\partial f_i}{\partial x_j} + \frac{\partial f}{\partial x_j} $.
By analyzing this we get the following result (this and similar results on forcing algebras can be found in \cite{gomezramirezthesis}).

\begin{proposition}
\label{singular}
Let $K$ be an algebraically closed field and let \[R=K[x_1, \ldots , x_m]/(g_1 , \ldots , g_k)\] be a $K$-algebra of finite type of dimension $d$. Let $f_1 , \ldots , f_n,f  \in R$ and let
\[B=R[t_1, \ldots, t_n]/(f_1t_1 + \ldots + f_nt_n+f) \]
be the forcing algebra. Let $P \in X=\Spec R$ be a closed point a,d let $Q=(P,t_1, \ldots, t_n)$ be a point over $P$. Then the following hold (in (3) and (4) we assume that $R$ and $B$ are domains, and that the forcing equation is not $0$).
\begin{enumerate}
\item
If there exists $f_i$ with $f_i(P) \neq 0$, then $Q$ is nonsingular if and only if $P$ is nonsingular.

\item
If $f_i(P)=0$ for all $i$ and $f(P) \neq 0$, then the fibre over $P$ is empty.

\item
If $f_i(P)=0$ for all $i$ and $f(P) = 0$, and if $P$ is a singular point, then $Q$ is also singular.

\item
If $f_i(P)=0$ for all $i$ and $f(P) = 0$, then every solution to the inhomogeneous linear system
\[  \begin{pmatrix}  \frac{\partial f_1}{\partial x_1}(P) & \ldots  & \frac{\partial f_n}{\partial x_1}(P) \\ \vdots & \ddots & \vdots  \\ \frac{\partial f_1}{\partial x_m}(P) & \ldots & \frac{\partial f_n}{\partial x_m}(P) \end{pmatrix}  \begin{pmatrix} t_1 \\ \vdots  \\  t_m \end{pmatrix}   = \begin{pmatrix} \frac{\partial f}{\partial x_1}(P)  \\ \vdots \\ \frac{\partial f}{\partial x_m}(P) \\ \end{pmatrix} \, \]
yields a singular point in the fiber over $P$.
\end{enumerate}
\end{proposition}
\begin{proof}
Statement (1) is clear, since the forcing algebra induces on each $D(f_i)$ an $(n-1)$-dimensional affine space. (2) is also clear. (3) and (4) follow from the Jacobian criterion for singularity \cite[Theorem I. 5.1]{hartshornealgebraic} (saying that a point $Q$ on an affine variety embedded in some affine space is nonsingular if and only if the rank of the Jacobian matrix $J(Q)$ equals the codimension of the variety). To prove (3), suppose that the rank of $J(P)$ is $< m - {\rm dim} \, R$. Then the rank of $J(Q)$ is at most $ m- {\rm dim} \, R < m+n -( {\rm dim }\, R +n-1)  =  m+n - {\rm dim }\, B$, hence $Q$ is a singular point. To prove (4), just observe that under this condition the rank $r$ of $J(P)$ and $J(Q)$ are the same. Hence the statement follows from $r \leq m -{\rm dim }\, R < m+n - ({\rm dim} \,R +n-1 )$.
\end{proof}

\begin{example}
Let $R=K[x,y,z]/(x+x^2+y^2+z^2)$ and $f_1=x^2$, $f_2=y^2$ and $f=x$. Then for $P=(0,0,0)$ neither condition (3) nor (4) of Proposition \ref{singular} is fulfilled. However, the rank of the Jacobian matrix $J(Q)$ for every Point $Q=(0,0,0,u,v)$ over $P$ is $1$ and therefore these points are all singular.
\end{example}

\bibliographystyle{amsplain}
\bibliography{bibliothek}

\providecommand{\bysame}{\leavevmode\hbox to3em{\hrulefill}\thinspace}
\providecommand{\MR}{\relax\ifhmode\unskip\space\fi MR }
\providecommand{\MRhref}[2]{%
  \href{http://www.ams.org/mathscinet-getitem?mr=#1}{#2}
}
\providecommand{\href}[2]{#2}
\begin{thebibliography}{10}

\bibitem{abetotalspace2}
M.~Abe, \emph{Steinness of the total space of a non-trivial algebraic affine
  {$\CC$}-bundle on the punctured complex affine plane}, Math. Nachr.
  \textbf{238} (2002), 16--22.

\bibitem{bingenerstorchdivisor}
J.~Bingener and U.~Storch, \emph{Zur {B}erechnung der {D}ivisorenklassengruppen
  kompletter lokaler {R}inge}, Nova acta Leopoldina \textbf{NF 52} (1981),
  no.~240, 7--63.

\bibitem{bonnetsurjectivity}
P.~Bonnet, \emph{Surjectivity of quotient maps for algebraic {$({\mathbb
  C},+)$}-actions and polyomial maps with contractible fibers}, Transformations
  Groups \textbf{7} (2002), no.~1, 3--14.

\bibitem{brennerdissertation}
H.~Brenner, \emph{{Zur Affinit\"at von Hyperfl\"achenkomplementen in affinen
  und projektiven Schemata}}, Dissertation, Ruhr-Universit\"at Bochum, 1998.

\bibitem{brennersuperheight}
\bysame, \emph{On superheight conditions for the affineness of open subsets},
  J. Algebra \textbf{247} (2002), 37--56.

\bibitem{brennerparasolid}
\bysame, \emph{How to rescue solid closure}, J. Algebra \textbf{265} (2003),
  579--605.

\bibitem{brennertightproj}
\bysame, \emph{Tight closure and projective bundles}, J. Algebra \textbf{265}
  (2003), 45--78.

\bibitem{brennerbarcelona}
\bysame, \emph{Tight closure and vector bundles}, Three Lectures on Commutative
  Algebra (J.~Elias S.~Zarzuela G.~Colom\'{e}-Nin, T. Cortodellas~Benitez,
  ed.), University Lecture Series, vol.~42, AMS, 2008, pp.~1--71.

\bibitem{brennerforcingsyzygy}
\bysame, \emph{Forcing algebras, syzygy bundles, and tight closure},
  Commutative Algebra. Noetherian and non-Noetherian Perspectives, Springer,
  2011.

\bibitem{brennerkolkata}
\bysame, \emph{Vector bundles and their torsors ({Kolkata} 2011)}, Lecture
  notes available on de.wikiversity.org (2011).

\bibitem{brunsherzog}
W.~Bruns and J.~Herzog, \emph{Cohen-{M}acaulay rings, revised edition},
  Cambridge University Press, 1998.

\bibitem{duboulozfinstonexotic}
A.~Dubouloz and D.~R. Finston, \emph{On exotic affine $3$-spheres}, Preprint
  (2011).

\bibitem{duboulozfinstonmehta}
A.~Dubouloz, D.~R. Finston, and P.D. Mehta, \emph{Factorial threefolds with
  {$G_a$}-actions}, Preprint (2009).

\bibitem{duttaguptaonoda}
A.K. Dutta, N.~Gupta, and N.~Onoda, \emph{Some patching results on algebras
  over two-dimensional factorial domains}, Preprint (2010).

\bibitem{finstonmaubach}
D.~R. Finston and S.~Maubach, \emph{The automorphism group of certain factorial
  threefolds and a cancellation problem}, Israel J. Math. \textbf{163} (2008),
  no.~1, 369--381.

\bibitem{flennerrationalquasihomogen}
H.~Flenner, \emph{Rationale quasi-homogene {Singularit\"aten}}, Arch. Math.
  \textbf{36} (1981), 35--44.

\bibitem{gomezramirezthesis}
Danny {Gomez Ramirez}, \emph{Normality properties of forcing algebras}, Ph.D.
  thesis, 2011, In preparation.

\bibitem{EGAII}
A.~Grothendieck and J.~Dieudonn\'{e}, \emph{El\'{e}ments de g\'{e}om\'{e}trie
  alg\'{e}brique {II}}, vol.~8, Inst. Hautes \'{E}tudes Sci. Publ. Math., 1961.

\bibitem{gurjarmasudamiyanishi}
R.V. Gurjar, K.~Masuda, and M.~Miyanishi, \emph{{${\mathbb A}^1$}-fibrations on
  affine threefolds}, Preprint (2010).

\bibitem{hararationalfrobenius}
N.~Hara, \emph{A characterization of rational singularities in terms of
  injectivity of {F}robenius maps}, Amer. J. of Math. \textbf{120} (1998),
  no.~5, 981--996.

\bibitem{hartshornealgebraic}
R.~Hartshorne, \emph{Algebraic {G}eometry}, Springer, New York, 1977.

\bibitem{hochstersolid}
M.~Hochster, \emph{Solid closure}, Contemp. Math. \textbf{159} (1994),
  103--172.

\bibitem{hochsterhunekebriancon}
M.~Hochster and C.~Huneke, \emph{Tight closure, invariant theory, and the
  {B}rian\c{c}on-{S}koda theorem}, J. Amer. Math. Soc. \textbf{3} (1990),
  31--116.

\bibitem{hunekeapplication}
C.~Huneke, \emph{Tight closure and its applications}, CBMS Lecture Notes in
  Mathematics, vol.~88, AMS, Providence, 1996.

\bibitem{lipmanrational}
J.~Lipman, \emph{Rational singularities}, Pub. Math. I.H.E.S. \textbf{36}
  (1969), 195--279.

\bibitem{samuelfactorial}
P.~Samuel, \emph{Unique factorization}, Amer. Math. Monthly \textbf{75} (1968),
  945--952.

\bibitem{schenzelflatfinite}
P.~Schenzel, \emph{When is a flat algebra of finite type}, Proc. A.M.S.
  \textbf{109} (1990), 287--290.

\bibitem{smithrational}
K.~E. Smith, \emph{F-rational rings have rational singularities}, Amer. J.
  Math. \textbf{119} (1997), 159--180.

\bibitem{watanaberational}
K.I. Watanabe, \emph{Rational singularities with {$k$}-action}, Lect. Notes in
  Pure and Appl. Math., vol.~84, 1983, pp.~339--351.

\bibitem{winkelmannfreeholomorphic}
J.~Winkelmann, \emph{On free holomorphic {${\mathbb C}$}-actions on {${\mathbb
  C}^n$} and homogeneous {S}tein manifolds}, Math. Ann. \textbf{286} (1990),
  593--612.

\end{thebibliography}

\end{document}